\documentclass[11pt]{amsart}

\usepackage{amsmath,amsthm,amssymb,amscd}

\newcommand{\p}{{\mathfrak{p}}}

\renewcommand{\r}{{\mathfrak{r}}}

\newcommand{\m}{{\mathfrak{m}}}

\newcommand{\Rbar}{\widetilde{R}}

\newcommand{\Spec}{\operatorname{Spec}}
\newcommand{\End}{\operatorname{End{}}}
\newcommand{\Hom}{\operatorname{Hom{}}}
\newcommand{\Ext}{\operatorname{Ext{}}}

\newcommand{\add}{\operatorname{add{}}}

\newcommand{\depth}{\operatorname{depth}}

\newcommand{\syz}{\operatorname{syz{}}}

\renewcommand{\phi}{\varphi}
\renewcommand{\to}{{\longrightarrow}}
\renewcommand{\mod}{\operatorname{mod}}

\newcommand{\coker}{\operatorname{coker}}

\newcommand{\gldim}{\operatorname{gldim}}
\newcommand{\repdim}{\operatorname{repdim}}
\renewcommand{\mod}{\operatorname{mod}}
\newcommand{\E}{\mathcal{E}}
\newcommand{\A}{\mathcal{A}}
\newcommand{\e}{\operatorname{e}}
\newcommand{\C}{\mathbb{C}}

\newcommand{\GL}{\operatorname{GL}}

\theoremstyle{plain}
\newtheorem{thm}{Theorem}
\newtheorem{cor}[thm]{Corollary}
\newtheorem*{cor*}{Corollary}
\newtheorem{prop}[thm]{Proposition}
\newtheorem{lemma}[thm]{Lemma}

\newtheorem*{AusA}{Theorem A}
\newtheorem*{AusB}{Theorem B}

\theoremstyle{remark}
\newtheorem{rem}[thm]{Remark}
\newtheorem{eg}[thm]{Example}

\theoremstyle{definition}
\newtheorem{defn}[thm]{Definition}

\begin{document}

\title{Endomorphism rings of finite global dimension}

\author{Graham J. Leuschke}
\address{Mathematics Department\\
	Syracuse University \\
	Syracuse NY 13244}
\email{gjleusch@syr.edu}
\urladdr{http://www.leuschke.org/}

\date{25 June 2004}

\thanks{This research began at MSRI while the author was participating
  in the 2002-3 program on Commutative Algebra.  Thanks to
  S.~Paul~Smith for comments on an earlier draft, and particular
  gratitude to Ragnar-Olaf Buchweitz for encouragement and suggestions
  over the course of several conversations on this topic.}
\bibliographystyle{amsplain}


\begin{abstract}  
For a commutative local ring $R$, consider (noncommutative)
$R$-algebras $\Lambda$ of the form $\Lambda = \operatorname{End}_R(M)$
where $M$ is a reflexive $R$-module with nonzero free direct summand.
Such algebras $\Lambda$ of finite global dimension can be viewed as
potential substitutes for, or analogues of, a resolution of
singularities of $\operatorname{Spec} R$.  For example, Van den Bergh
has shown that a three-dimensional Gorenstein normal
$\mathbb{C}$-algebra with isolated terminal singularities has a
crepant resolution of singularities if and only if it has such an
algebra $\Lambda$ with finite global dimension and which is maximal
Cohen--Macaulay over $R$ (a ``noncommutative crepant resolution of
singularities'').  We produce algebras
$\Lambda=\operatorname{End}_R(M)$ having finite global dimension in
two contexts: when $R$ is a reduced one-dimensional complete local
ring, or when $R$ is a Cohen--Macaulay local ring of finite
Cohen--Macaulay type.  If in the latter case $R$ is Gorenstein, then
the construction gives a noncommutative crepant resolution of
singularities in the sense of Van den Bergh.
\end{abstract}

\maketitle
\numberwithin{thm}{section}
\numberwithin{equation}{thm}

This paper takes for its starting point two results of Auslander:

\begin{AusA}[{\cite[\S III.3]{Auslander:QueenMary}}] 
Let $\Lambda$ be a left Artinian ring with radical $\r$ and assume
that $\r^n=0$, $\r^{n-1}\neq 0$.  Set $M = \bigoplus_{i=1}^n
\Lambda/\r^i$.  Then $\Gamma := \End_\Lambda(M)^{op}$ is coherent of
global dimension at most $n+1$. \end{AusA}

\begin{AusB}[{\cite{Auslander:rationalsing}}] 
Let $S=k[\![x,y]\!]$ be the ring of formal power series in two
variables over a field $k$ and let $G$ be a finite subgroup of
$\GL_2(k)$ with $|G|$ invertible in $k$.  Set $R=S^G$.  Then
$A:=\End_R(S)^{op}$ has global dimension at most two.  \end{AusB}

These theorems both relate to Auslander's notion of {\it
representation dimension}, introduced in \cite{Auslander:QueenMary} as
a way to measure homologically the failure of an Artin algebra to have
finite representation type.  The representation dimension of an Artin
algebra $\Lambda$ can be defined as
$$\repdim\Lambda = \inf \{ \gldim \End_\Lambda(M)\},$$  where the
infimum is taken over all finitely generated modules $M$ which are
generator-cogenerators for $\mod\Lambda$.  Note that Theorem~A does
not prove finiteness of the representation dimension; while $M$ has a
nonzero free direct summand, it need not be a cogenerator.  Auslander
showed in \cite{Auslander:QueenMary} that $\repdim\Lambda \leq 2$ if
and only if $\Lambda$ has finite representation type, and in 2003,
Rouquier \cite{Rouquier:2003} gave the first examples with
$\repdim\Lambda > 3$.  Iyama has recently shown \cite{Iyama:finrepdim}
that the representation dimension of an Artin algebra is always finite.

We extend Theorems~A and B in two directions.  In each case, we 
consider commutative Noetherian (semi)local base rings.

First, we fill a gap between Auslander's theorems: the case of
dimension one.  A reduced complete local ring $R$ of dimension one
always has a finitely generated module whose endomorphism ring has
finite global dimension: the normalization $\Rbar$.  However, $\Rbar$
is never a generator in the category of $R$-modules, unless $R$ is
already a discrete valuation ring.
Theorem~\ref{dimone} produces a finitely generated generator $M$ such 
that $\End_R(M)$ has finite global dimension.  Specifically, 
$M$ can be taken to be a direct sum of certain overrings $S$ between 
$R$ and $\Rbar$, and the global dimension of $\End_R(M)$ is bounded 
by the multiplicity of $R$.  This completes a coherent picture for 
rings of dimension at most 2; see \cite{VandenBergh:crepant} and 
\cite{VandenBergh:flops} for related progress in dimension three.

We also generalize to dimensions $d > 2$ by exploiting the connection
with finite representation type.  The two-dimensional quotient
singularities $\C[\![x,y]\!]^G$, with $G \subseteq \GL_2(\C)$ a finite
group, are precisely the two-dimensional complete local rings with
residue field $\C$ and having finite Cohen--Macaulay type
\cite{Herzog:1978, Auslander:rationalsing, Esnault:1985}.  Moreover,
$\C[\![x,y]\!]$ contains as $R$-direct summands all indecomposable
maximal Cohen--Macaulay $R$-modules \cite{Herzog:1978}.
Theorem~\ref{fcmt} states that if $R$ is a $d$-dimensional
Cohen--Macaulay local ring of finite Cohen--Macaulay type, and if $M$
is the direct sum of all indecomposable maximal Cohen--Macaulay
$R$-modules, then $\End_R(M)$ has global dimension at most $\max\{2,
d\}$.  It follows that the representation dimension of a complete
Cohen--Macaulay local ring of finite Cohen--Macaulay type is finite.
(See Definition~\ref{defrepdim} for the definition of representation
dimension in this context.)

The proofs of both theorems are based on projectivization
\cite[II.2]{AusReitSmalo}.  In the present contexts,  this means that
the functor $\Hom_R(M,-)$ induces an equivalence of categories between
$\add(M)$, the full subcategory of $R$-modules which are direct
summands of finite direct sums of copies of $M$, and the full
subcategory of finitely generated projective modules over $A
:=\End_R(M)^{op}$.  In particular, if $M$ is an $R$-generator, then
$M$ has an $R$-free direct summand, and so $M$ is a projective
$A$-module.  Auslander's original proof of Theorem~A uses this
technique, and the proofs of our two main results are very close in
spirit to his method.  Iyama has refined Auslander's methods into a
theory of {\it rejective subcategories} \cite{Iyama:rejective, 
Iyama:repdim_Solomon}, which he uses to prove that the representation
dimension of an order over a complete discrete valuation ring is
finite.

Section~\ref{connections} discusses the implications of
Theorem~\ref{fcmt} to the theory of non-commutative crepant
resolutions \cite{VandenBergh:flops}.  If $R$ is Gorenstein and has
finite Cohen--Macaulay type, Theorem~\ref{fcmt} does indeed produce a
non-commutative crepant resolution of $R$.  If $R$ is not Gorenstein,
then non-commutative crepant resolutions are not yet defined, but
Theorem~\ref{fcmt} still gives an analogue.  We discuss advantages and
disadvantages of this analogy.

The rings under consideration will be Noetherian, and all modules 
finitely generated. We abbreviate $\Hom_R(-,-)$ by $(-,-)$.

\section{Dimension one}

In this section we consider reduced one-dimensional semilocal rings
$R$.  We always assume that $R$ is complete with respect to its
Jacobson radical, equivalently
that $R$ is isomorphic to a direct product of complete local rings.
Let $K$ be the total quotient ring of $R$, obtained by inverting all
nonzerodivisors of $R$.  Recall that a finitely generated $R$-module
$M$ is {\it torsion-free} provided the natural map $M \to M\otimes_R
K$ is injective.

\medskip
Our goal requires us to consider the module theory of certain
birational extensions of reduced rings, that is, extensions $R
\subseteq S$ where $S$ is a finitely generated $R$-module contained in
the total quotient ring $K$ of $R$.  Of course, in this situation
every finitely generated torsion-free $S$-module is a finitely
generated torsion-free $R$-module, but not vice versa.  The following
lemma, however, follows easily by clearing denominators.

\begin{lemma}\label{homs}  Let $R \subseteq S$ be a birational extension 
of reduced rings as above. Let $C$ and $D$ be finitely generated
torsion-free $S$-modules.  Then $\Hom_R(C,D)=\Hom_S(C,D)$.
Furthermore, if $M$ is a finitely generated torsion-free $R$-module,
and $f :C \to M$ is an $R$-linear map, then the image of $f$ is an
$S$-module.\end{lemma}

For the remainder of this section, $(R, \m)$ will be a reduced
complete local ring of dimension one with total quotient ring $K$ and
integral closure $\Rbar$.  Note that $K$ is a direct product of
finitely many fields, and $\Rbar$ is correspondingly a direct product
of discrete valuation rings.  Since $R$ is complete and reduced,
$\Rbar$ is a finitely generated $R$-module \cite[Theorem
11.7]{Matsumura}.

Set $R^{(1)} := \End_R(\m)$.  Since $\m$ contains a nonzerodivisor,
$R^{(1)}$ embeds naturally into $K$ (by sending $f$ to $f(r)/r$, which
is independent of the nonzerodivisor $r$).  It is well known that in
fact $R^{(1)} \subseteq \Rbar$.  Furthermore, $R \subsetneq R^{(1)}$
unless $R = \Rbar$.  Now, $R^{(1)}$ may no longer be local (if, for
example, we take $R = k[\![x,y]\!]/(xy)$), but by Hensel's Lemma,
$R^{(1)}$ is a direct product of complete local rings, $R^{(1)} =
 R^{(1)}_1 \times \dots \times R^{(1)}_{n_1}$, each of which is again 
reduced.

Iterating this procedure by taking the endomorphism ring of the
maximal ideal of each of the local rings $R^{(1)}_\ell$, $\ell = 1,
\dots, n_1$, gives a family of reduced complete local rings
$\{R^{(i)}_j\}$.  Since $\Rbar/R$ is an $R$-module of finite length
and each $R_j^{(i)}$ is trapped between some $R_l^{(i-1)}$ and
$\Rbar$, this family is finite.  It follows that the lengths of the
chains
\begin{equation}\label{ringchain}\tag{$\ddagger$}
R \subsetneq R^{(1)}_{j_1} \subsetneq \dots \subsetneq R^{(n)}_{j_n} = \Rbar,
\end{equation}  
are bounded above, where each $R^{(i)}_{j_i}$ is a direct factor of 
the endomorphism ring of the maximal ideal of $R^{(i-1)}_{(j_{i-1})}$.

Let $\E(R)$ denote the family of rings obtained in this way, including
$R$ itself.   Put $\A(R) = \add(\E(R))$, the full subcategory of
$\operatorname{mod-}\!R$ containing all direct summands of finite
direct sums of rings in $\E(R)$, considered as $R$-modules.  If $S$ is
a direct product of complete local rings $S_i$, $j=1, \dots, m$, let
$\E(S)$ be the corresponding union of the $\E(S_j)$, and $\A(S) =
\add(\E(S))$ the full subcategory containing all direct summands of
finite direct sums of rings in $\E(S)$, again considered as
$S$-modules.

Even though we begin with a local ring, the proof of
Theorem~\ref{dimone} requires dealing with semilocal rings that crop
up along the way.  Lemma~\ref{semilocal} allows us to reduce to the
local case each time.

\begin{lemma}\label{semilocal} Let $S = S_1 \times \dots \times S_k$ 
be a direct product of rings.  Assume that for each $i = 1, \dots, k$ 
and for each torsion-free $S_i$-module $D$, there is an exact sequence  
\begin{equation}\label{semilocal1}
0 \to C_{i,m_i} \to C_{i,{m_i}-1} \to \cdots \to C_{i,0} \to D \to 0
\end{equation}
with each $C_{ij} \in \A(S_i)$ and such that 
\begin{equation*}
0 \to (X, C_{i,m_i}) \to (X, C_{i,{m_i}-1}) \to \cdots \to 
(X, C_{i,0}) \to (X, D) \to 0
\end{equation*}
is exact for all $X \in \A(S_i)$. Then for each torsion-free $S$-module 
$E$, there exists an exact sequence
\begin{equation}\label{semilocal3}
0 \to C_m \to C_{m-1} \to \cdots \to C_0 \to E \to 0
\end{equation}
with each $C_j \in \A(S_i)$ and such that 
\begin{equation}\label{semilocal4}
0 \to (X, C_m) \to (X, C_{m-1}) \to \cdots \to (X, C_0) \to (X, E) \to 0
\end{equation}
is exact for all $X \in \A(S)$.
\end{lemma}

\begin{proof} Let $E$ be a torsion-free $S$-module.  Then $E \cong 
\prod_{i=1}^n e_iE$, where $e_i$ is a complete set of orthogonal 
idempotents for the decomposition $S = S_1 \times \dots \times S_n$.  
The exact sequence (\ref{semilocal3}) can be taken to be the direct 
sum of the sequences (\ref{semilocal1}) with $D=e_i E$.  It remains 
to show that (\ref{semilocal4}) is exact for all $X \in \A(S)$.  Since 
$\Hom_S(Y, Z)=0$ whenever $Y$ is a $S_i$-module and $Z$ is a $S_j$-
module with $i \neq j$, this is clear.
\end{proof}

We can now state the key result which will imply our main theorem in 
the reduced case. 

\begin{prop}\label{keyred} Let $(R, \m)$ be a reduced complete local 
ring of dimension one and let $N$ be a torsion-free $R$-module.  Let 
$n$ be the length of the longest chain (\ref{ringchain}). Then there 
exists an exact sequence
\begin{equation*}
0 \to C_n \to C_{n-1} \to \cdots \to C_0 \to N \to 0
\end{equation*}
with each $C_i \in \A(R)$ and such that 
\begin{equation*}
0 \to (X, C_n) \to (X, C_{n-1}) \to \cdots \to (X, C_0) \to (X, N) \to 0
\end{equation*}
is exact for all $X \in \A(R)$.
\end{prop}

\begin{proof} We proceed by induction on $n$.
If $n=0$, then $R = \Rbar$ is a discrete valuation ring, and any 
torsion-free $R$-module is free.  The set $\A(R)$ consists exactly 
of the free $R$-modules, so that $0 \to C_0 \overset{=}\to N \to 0$ 
is the required sequence.

Assume that the statement holds for reduced complete local rings of 
dimension one having a chain (\ref{ringchain}) of length at most $n-1$, 
and that $R$ has a chain of length $n$.  In particular, then the 
proposition holds for each direct factor of $R^{(1)} = \End_R(\m)$.  
Let $N$ be a torsion-free $R$-module.

First suppose that 
$N$ is an $R^{(1)}$-module.  By 
Lemma~\ref{semilocal}, then, there is an exact sequence of $R^{(1)}$-modules
\begin{equation}
0 \to C_n \to C_{n-1} \to \cdots \to C_0 \to N \to 0,
\end{equation}
with each $C_i \in \A(R^{(1)})$, which remains exact under $(X, -)$
for any $X \in \A(R^{(1)})$.  (We use Lemma~\ref{homs} here to know
that $\Hom_R(X, -) = \Hom_{R^{(1)}}(X,-)$.)  The only indecomposable
module in $\A(R)$ but not in $\A(R^{(1)})$ is the free module $R$, so
the sequence remains exact under $(X, -)$ for any $X \in \A(R)$, as
desired.

Next suppose that $N$ is not an $R^{(1)}$-module.  Let $N'
=\Hom_R(R^{(1)}, N) \subsetneq N$ be the largest $R^{(1)}$-module
contained in $N$.  Observe that for any $R^{(1)}$-module $X$, and any
$R$-linear homomorphism $X \to N$, the image of $X$ is contained in
$N'$.  In particular, $(X, N') = (X,N)$ for any $X \in \A(R^{(1)})$.
By induction, there is a surjection $f: C' \to N'$, with $C' \in
\A(R^{(1)})$, such that $(X, f): (X, C') \to (X, N')$ is surjective
for all $X \in \A(R^{(1)})$.  Since $(X, N') = (X,N)$, we see that
applying $(X,-)$ to the composition $C' \to N' \hookrightarrow N$
yields a surjection for all $X \in \A(R^{(1)})$.

Take a free $R$-module $F$ mapping minimally onto $N/N'$ and lift to a
homomorphism $g: F \to N$.  Then $g^{-1}(N')$ is an $R^{(1)}$-module.
Indeed, $g^{-1}(N') \subseteq \m F$ as $F$ is a minimal free cover of
$N/N'$, and since $\m F$ is clearly a module over
$R^{(1)}=\End_R(\m)$, Lemma~\ref{homs} implies that $g^{-1}(N') =\m
F$, so in particular is an $R^{(1)}$-module.

Define $\pi: F \oplus C' \to N$ by $\pi(p,c) = g(p) - f(c)$.  Since
$(X, f)$ is surjective for all $X \in \A(R^{(1)})$, and $g$ induces a
surjection $F \to N/N'$, we see that $(X, \pi)$ is surjective for all
$X \in \A(R)$.  We claim that the $L =\ker \pi$ is an
$R^{(1)}$-module.   Let $\alpha \in R^{(1)}$ and $(p,c) \in L$, so
that $g(p) = f(c)$.  Since $f(c) \in N'$ and $g^{-1}(N')$ is an
$R^{(1)}$-module, $\alpha p \in g^{-1}(N')$.  Then, since
$f|_{f{-1}(N')}$ is $R^{(1)}$-linear by Lemma~\ref{homs}, $g(\alpha p)
= \alpha g(p) =\alpha f(c)$.  Finally, $f: C' \to N'$ is
$R^{(1)}$-linear, so $\alpha f(c) = f(\alpha c)$.  That is, $(\alpha
p, \alpha c) \in L$, as claimed.

By the previous case, then, there is an exact sequence
\begin{equation*}
0 \to C_{n-1} \to C_{n-2} \to \cdots \to C_0 \to L \to 0
\end{equation*}
such that 
\begin{equation*}
0 \to (X, C_{n-1}) \to (X, C_{n-2}) \to \cdots \to (X, C_0) \to (X,L)\to 0
\end{equation*}
is exact for all $X \in \A(R^{(1)})$. Splicing this together with the
short exact sequence $0 \to L \to F \oplus C' \to N \to 0$, and using
once again that the only indecomposable module in $\A(R)$ but not in
$\A(R^{(1)})$ is the free module $R$, we are done.
\end{proof}

\begin{thm}\label{dimone} Let $(R,\m)$ be a one-dimensional reduced 
complete local ring. Put $M =\bigoplus_{S \in \E(R)} S$, a finitely
generated $R$-module.  Then $\Gamma := \End_R(M)^{op}$ has global
dimension at most $n+1$, where $n$ is the length of the longest chain
(\ref{ringchain}).\end{thm}

\begin{proof} Let $N$ be a finitely generated $\Gamma$-module. Then by 
\cite[II.2]{AusReitSmalo} there exists a homomorphism $M_1 \to M_0$
with $M_i \in\add(M)$ such that $(M, M_1) \to (M, M_0) \to N \to 0$ is
exact.  Let $L$ be the kernel of  $M_1 \to M_0$.  Then $L$ is
torsion-free, so by the Proposition has a resolution of length $n$ by
modules in $\add(M)$, which remains exact after applying $(M,-)$.
Since each $(M,M_i)$ is a projective $\Lambda$-module, $N$ has
projective dimension at most $n+1$.
\end{proof}

\begin{cor} A one-dimensional reduced complete local ring has a finitely 
generated module whose endomorphism ring has global dimension at most
$\e(R)$, the multiplicity of $R$.
\end{cor}

\begin{proof} It is known that $\Rbar/\m\Rbar$ has dimension $\e(R)$ 
as a vector space over $R/\m$.  Thus $\Rbar/R$ has length $\e(R)-1$, 
so $e$ is a uniform bound on the length of chains (\ref{ringchain}).
\end{proof}

It seems plausible that Theorem~\ref{dimone} actually holds for rings
$R$ such that the integral closure $\Rbar$ is a finitely-generated
$R$-module and a regular ring, for example, the hypersurface
$x^2+y^3-y^2z^2 = 0$.  The proof given above is reminiscent of the
algorithm of de Jong \cite{deJong} for obtaining the integral closure
by taking iterated endomorphism rings.

\section{Finite Cohen--Macaulay Type}

As mentioned in the introduction, the original motivation for
Auslander's representation dimension was to study Artin algebras of
finite representation type, that is, Artin algebras with only finitely
many isomorphism classes of finitely generated modules.  For
(commutative Noetherian) rings of higher dimension, this property has
been generalized to {\it finite Cohen--Macaulay type}.  A nonzero
finitely generated module $M$ over a $d$-dimensional ring $R$ is
called maximal Cohen--Macaulay (MCM) if there exists an $M$-regular
sequence $x_1, \ldots, x_d$.  We then say that $R$ has finite CM type
provided there are, up to isomorphism, only finitely many
indecomposable MCM $R$-modules.

The one-dimensional CM local rings of finite CM type are completely
characterized \cite{Drozd-Roiter, Green-Reiner, Wiegand:1989,
Wiegand:1994, Cimen:thesis}.  In dimension two, the complete local
rings containing the complex numbers and having finite CM type are
also completely classified \cite{Herzog:1978, Auslander:rationalsing,
Esnault:1985}.  They are exactly the rings of Theorem~B, that is, the
invariant rings $R = \C[\![x,y]\!]^G$ under the action of a finite
group $G$. In this case, $\C[\![x,y]\!]$ is a representation generator
for $R$, that is, contains as direct summands all the indecomposable
MCM $R$-modules.

The main result of this section is a generalization of Theorem~B.
Again, the proof relies on the process of projectivization, which was
described in the introduction.

\begin{thm}\label{fcmt} Let $(R, \m)$ be a $d$-dimensional CM local 
ring of finite CM type.  Let $M$ be a representation generator for 
$R$ (in particular, $M$ has a free direct summand).  Then $A := 
\End_R(M)^{op}$ has global dimension at most $\max\{2, d\}$.  If $d 
\geq 2$, then equality holds.
\end{thm}

\begin{proof} First assume that $d\geq 2$.  Let $N$ be a finitely 
generated left $A$-module.  Take the first $d-1$ steps in a projective 
resolution of N over A:
\begin{equation}
\CD \textbf{P}_\bullet: P_{d-1} @>\phi_{d-1}>> \cdots @>\phi_2>> P_1
@>\phi_1>> P_0 \endCD
\end{equation}
with $\coker \phi_1 = N$.   For each $i$, we can use
\cite[II.2]{AusReitSmalo} to write $P_i =\Hom_R(M,M_i)$, where $M_i$
is a direct summand of a direct sum of copies of $M$. In particular,
each $M_i$ is a MCM $R$-module.  Moreover, each $\phi_i$ can be
written as $\Hom_R(M,f_i)$ for $R$-homomorphisms $f_i: M_i \to
M_{i-1}$.  This gives the following sequence of MCM $R$-modules and
homomorphisms:
\begin{equation}
\CD
\textbf{C}_\bullet: M_{d-1} @>f_{d-1}>> \cdots @>f_2>> M_1 @>f_1>> M_0.
\endCD
\end{equation}
Since $M$ has a free direct summand and $\Hom_R(M,
\textbf{C}_\bullet)$ is exact, it follows that in fact
$\textbf{C}_\bullet$ is exact.  Put $M_d = \ker(f_{d-1})$.  Then $M_d$
is a MCM $R$-module by the depth lemma, and left-exactness of $\Hom$
gives an exact sequence
\begin{equation}
\CD 0 @>>> \Hom_R(M,M_d) @>>> P_{d-1} @>\phi_{d-1}>> \cdots @>\phi_2>>
P_1 @>\phi_1>> P_0.  \endCD
\end{equation}
Since $M_d$ is MCM, $\Hom_R(M,M_d)$ is $A$-projective, and $N$ has
projective dimension at most $d$.

To see that the global dimension of $A$ is exactly $d$, take $N$ to be
a simple $A$-module.  Then $N$ has finite length as an $R$-module, so
is of depth zero.  A projective resolution of $N$ is in particular an
exact sequence of MCM $R$-modules, and so a projective resolution of
length less than $d$ would contradict the depth lemma.  Thus $N$ has
projective dimension exactly $d$.  Since the global dimension of $A$
is the maximum of the projective dimensions of the simple modules,
this finishes the case $d\geq 2$.

If $d<2$, then we can repeat the first part of the argument, simply
taking an $A$-projective resolution of length 1.  The remainder of the
proof is the same, showing that $\gldim A \leq 2$.
\end{proof}

One can also state the proof above in terms of the adjoint pair $g_* =
\Hom_R(M,-)$ and $g^*=-\otimes_A M$.  In this context, the fact that
$\textbf{C}_\bullet$ is exact comes down to the facts that (1) $M$ is
a projective $A$-module, so $g^*$ is an exact functor, and (2)
$g_*g^*$ is the identity on projective $A$-modules.

Following Iyama \cite{Iyama:repdim_Solomon}, we extend the definition
of representation dimension to rings of positive Krull dimension.

\begin{defn}\label{defrepdim} Let $T$ be a complete regular local ring 
and let $R$ be a $T$-algebra, finitely generated and free as a 
$T$-module.  Let $\mathcal C$ be the category of $R$-modules which 
are finitely generated free $T$-modules.  Let 
$$\repdim R = \underset{N \in \mathcal C}{\inf} \{ \gldim \End_R(R
\oplus R^* \oplus N)\},$$ where $R^* := \Hom_T(R,T)$.\end{defn}

\begin{prop} Let $R$ be a CM complete local ring of finite CM type.  
Then $\repdim R \leq \max\{2, \dim R\}$.\end{prop}

\begin{proof} By Cohen's structure theorem, $R$ is a finitely generated 
module over some complete regular local ring $T$.  The MCM $R$-modules 
are precisely the $R$-modules that are free over $T$.  Finally, $R^* 
\cong \omega_R$ is the canonical module for $R$, which is MCM.  
Theorem~\ref{fcmt} then shows that $\End_R(R \oplus \omega_R \oplus N)$ 
has global dimension at most $\max\{2, \dim R\}$, where $N$ is the 
direct sum of the remaining indecomposable MCM $R$-modules.
\end{proof}

\begin{rem} 
As mentioned in the Introduction, Auslander proved in
\cite{Auslander:QueenMary} that if $\Lambda$ is an Artin algebra of
finite representation type, with additive generator $M$, then $\Gamma
:= \End_\Lambda(M)$ has representation dimension two.  In fact, he
showed that $\Gamma$ is what is now called an {\it Auslander
  algebra\/}, that is, $\Gamma$ has global dimension two and dominant
dimension two.  We say that a ring $A$ has {\it dominant dimension} at
least $t$ if, in a minimal injective resolution $0 \to A \to
I^\bullet$ of $A$, the $j^\text{th}$ injective module $I^j$ is also
projective for $j<t$. 

One checks easily that the proof given in \cite[VI.5]{AusReitSmalo}
applies verbatim to CM local rings of finite CM type, and shows that,
in the situation of Theorem~\ref{fcmt}, $A = \End_R(M)$ has dominant
dimension two.  To see that one cannot hope for higher dominant
dimension, consider a three-dimensional CM local ring $R$ of finite CM
type, {\it e.g.} Example~\ref{eg1} or \ref{eg2}.  Let $M$ be the
direct sum of the indecomposable MCM $R$-modules, and $A=\End_R(M)$.
If $A$ has dominant dimension $>2$, then the minimal injective
resolution of $A$ is 
$$
0 \to A \to \Hom_R(M,I^0) \to \Hom_R(M,I^1) \to \Hom_R(M,I^2)
$$
where each $\Hom_R(M,I^j)$ is a projective-injective $A$-module.  It
follows that each $I^j$ is an injective $R$-module.  By the
equivalence of categories between $\add(M)$ and projective
$A$-modules, then, the minimal injective resolution of $M$ over $R$ is 
$$
0 \to M \to I^0 \to I^1 \to I^2,
$$
and applying $\Hom_R(M,-)$ to this injective resolution preserves
exactness.  This implies that $\Ext_R^1(M,M)=0$, which is quite false
in Examples~\ref{eg1} and \ref{eg2}. 
\end{rem}

\section{Connections with non-commutative crepant resolutions}
\label{connections}

Though it is a purely algebraic statement, Theorem~\ref{fcmt} is
closely related to geometric statements about resolution of
singularities.  Recent work of M.~Van den Bergh
\cite{VandenBergh:crepant, VandenBergh:flops} has revealed unexpected
connections between the (geometric) resolutions of certain rational
singularities and the algebraic properties of certain endomorphism
rings over their coordinate rings.  To make this connection more
precise, we quote the following definition of Van den Bergh:

\begin{defn}\label{nccrepdef} Let $R$ be a Gorenstein normal domain.  
A {\em non-commutative crepant resolution} of $R$ is an $R$-algebra 
$A = \End_R(M)$, for a finitely generated reflexive $R$-module $M$, 
such that $A$ is a MCM $R$-module and $\gldim A_\p = \dim R_\p$ for 
all $\p \in \Spec R$.\end{defn}

Non-commutative crepant resolutions were introduced to solve a problem
related to their geometric counterparts.  A {\it (geometric) crepant
resolution} of a scheme $X$ is a projective morphism $f: Y \to X$,
with $Y$ regular, such that $f^*\omega_X = \omega_Y$, where $\omega$
denotes the canonical bundle.  A.~Bondal and D.~Orlov
\cite{Bondal-Orlov} conjecture that if $X$ has a geometric crepant
resolution, then any two have equivalent bounded derived categories of
coherent sheaves.  M.~Kapranov and E.~Vasserot
\cite{Kapranov-Vasserot} verify the conjecture of Bondal-Orlov for
two-dimensional quotient singularities $\C^2/G$; they show that any
geometric crepant resolution is derived equivalent to the
non-commutative crepant resolution given by $\End_R(\C[\![x,y]\!])$.

Van den Bergh pushes this point of view into dimension three.  He
proves \cite{VandenBergh:crepant} that a three-dimensional Gorenstein
normal $\C$-algebra with terminal singularities has a non-commutative
crepant resolution if and only if it has a geometric crepant
resolution, and furthermore that the two crepant resolutions are
derived equivalent, establishing the conjecture of Bondal--Orlov in
this case.  In fact, he conjectures
\cite[Conj. 4.6]{VandenBergh:crepant} that all crepant resolutions of
a given Gorenstein scheme $X$, non-commutative as well as geometric,
are derived equivalent.  He also gives some further examples in which
non-commutative crepant resolutions exist.

\nocite{VandenBergh:flops}

Theorem~\ref{fcmt} implies that if $R$ is a Gorenstein local ring of
finite CM type, containing a field and having dimension two or three,
then $R$ has a non-commutative crepant resolution.  The completion of
such a ring is the analytic local ring of one of the simple
hypersurface singularities (see, for example, \cite{Yoshino:book}),
the MCM modules of which are known.  One can thus check that for each
of the simple singularities, the endomorphism ring of a representation
generator is indeed a MCM module.

In fact, Theorem~\ref{fcmt} gives a little more: the reflexive module
in the definition of the non-commutative crepant resolution can
actually be taken to be MCM.  This is consistent with all the other
known examples of non-commutative crepant resolutions \cite[Remark
4.4]{VandenBergh:crepant}.

It is worth pointing out that Theorem~\ref{fcmt} also implies that
graded 
Gorenstein local rings of finite CM type (containing $\C$) have
rational singularities.  Van den Bergh shows
\cite[Prop. 3.3]{VandenBergh:crepant} that if $R$ is 
a positively graded Gorenstein algebra over a field, with
an isolated singularity, and $R$ has a noncommutative crepant resolution of
singularities, then $R$ has at most rational singularities.  Of
course, this is already known for the Gorenstein local rings of finite
CM type, following from the complete classification of such rings
\cite{Herzog:1978, BGS, Knorrer}.

Non-commutative crepant resolutions are not yet defined for
non-Gorenstein rings.  It is tempting to accept
Definition~\ref{nccrepdef} verbatim for all CM normal domains, and
look to Theorem~\ref{fcmt} as a source of non-commutative crepant
resolutions in this context as well.  This optimism is quickly
tempered by the following two examples (the only known non-Gorenstein
CM local rings of of finite CM type and dimension $\geq 3$).

\begin{eg}\label{eg1} 
Let $R = k[\![x^2, xy, y^2, yz, z^2, xz]\!]$, where $k$ is 
an algebraically closed field of characteristic zero.  By
\cite[16.10]{Yoshino:book}, $R$ has finite CM type.  The
indecomposable MCM $R$-modules are the free module of rank one, the
canonical module $\omega \cong (x^2, xy, xz)$, and
$M:=\syz_1^R(\omega)$, which has rank $2$.  By Theorem~\ref{fcmt},
$A:=\End(R\oplus\omega\oplus M)$ has global dimension $3$.  However,
$\depth_R A =2$ (this can be easily checked with, say, Macaulay2
\cite{M2}).  The culprit is $M$: both $\Hom_R(M,R)$ and $\Hom_R(M,M)$
have depth $2$.

Removing $M$, however, eliminates the problem. Observe that $R$ is a
ring of invariants of $k[\![x,y,z]\!]$ under an action of ${\mathbb
Z}_2$.  Therefore $\End_R(k[\![x,y,z]\!])$ is isomorphic to the
twisted group ring $k[\![x,y,z]\!] * {\mathbb Z}_2$, and the twisted
group ring has global dimension $3$ by \cite[Ch. 10]{Yoshino:book}.
Finally, since $k[\![x,y,z]\!] \cong R \oplus \omega$ as an
$R$-module, we see that $R\oplus \omega$ gives a noncommutative
crepant resolution of $R$, and exhibits $\repdim R \leq 3$. \end{eg}

\begin{eg}\label{eg2} 
Let $R = k[\![x,y,z,u,v]\!]/I$, where $I$ is generated by 
the $2\times 2$ minors of the matrix $\left(\smallmatrix x & y & u \\
y & z & v \endsmallmatrix\right)$.  Then $R$ has finite CM type
\cite[16.12]{Yoshino:book}.  The only indecomposable nonfree MCM
$R$-modules are, up to isomorphism,
\begin{itemize}
\item the canonical module $\omega \cong (u,v)R$;
\item $M := \syz_1^R(\omega)$, isomorphic  to the ideal $(x,y,u)R$;
\item $N:= \syz_2^R(\omega)$, rank two and $6$-generated;
\item $L:=M^\vee$, the canonical dual of $M$, isomorphic to the ideal
  $(x,y,z)R$. 
\end{itemize}
In particular, $\omega^* = \Hom_R(\omega,R)$ is isomorphic to $M$.  By
Theorem~\ref{fcmt}, then, $A := \End_R(R \oplus \omega \oplus M \oplus
N \oplus L)$ has global dimension $3$.  Again, $A$ fails to be MCM as
an $R$-module, since none of $L^*$, $N^*$, and $\Hom_R(\omega,M)$ are
MCM.  In this example, $\End_R(R\oplus\omega)$ and $\End_R(R\oplus M)$
are the only endomorphism rings of the form $\End_R(D)$, with $D$
nonfree MCM, that are themselves MCM.  I do not know whether the
endomorphism ring $\End_R(R\oplus\omega)$ has finite global dimension.
\end{eg}


\providecommand{\bysame}{\leavevmode\hbox to3em{\hrulefill}\thinspace}
\providecommand{\MR}{\relax\ifhmode\unskip\space\fi MR }
\providecommand{\MRhref}[2]{%
  \href{http://www.ams.org/mathscinet-getitem?mr=#1}{#2}
}
\providecommand{\href}[2]{#2}

\end{document}